\documentstyle[12pt]{amsart}

\makeatletter                                            
\renewcommand{\@begintheorem}[2]{                        
\rm \trivlist \item [\hskip \labelsep {\bf #2\ \ #1.}]   
                                }                        
\makeatother                                             

\makeatletter

\newcommand{\ts}{\vspace{\baselineskip}\noindent{\bf Proof.}$\;\;$}
\newcommand{\ZZ}{{\bf Z}}

\newcommand{\CC}{{\bf C}}

\newcommand{\PP}{{\bf P}}

\newcommand{\cO}{{\cal O}}
\newcommand{\cM}{{\cal M}}

\begin{document}

\title{Hessians and the moduli space of cubic surfaces}

\author{Elisa Dardanelli}
\address{Piazza Biancani 17, 12100 Cuneo, Italia}
 \email{eli.gian@tiscali.it}
\author{Bert van Geemen}
\address{Dipartimento di Matematica, Universit\`a di Milano,
  via Saldini 50, I-20133 Milano, Italia}
 \email{geemen@mat.unimi.it}

\begin{abstract}
The Hessian of a general cubic surface is a nodal quartic surface, hence
its desingularisation is a K3 surface.
We determine the transcendental lattice of the Hessian K3 surface
for various cubic surfaces (with nodes and/or Eckardt points for example).
Classical invariant theory shows that the moduli space of cubic surfaces
is a weighted projective space. We describe the singular locus
and some other subvarieties of the moduli space.

\end{abstract}

\maketitle

Moduli spaces of cubic surfaces have been intensively studied recently.
The main reason was the discovery by Allcock, Carlson and Toledo \cite{ACT}
that the moduli space of cubic surfaces is a ball quotient.
Another proof of this fact, using K3 surfaces rather than cubic
threefolds, is given in \cite{DGK}.
There is a classical, quite different, way to associate a K3 surface to
a general cubic surface,
it is simply by taking the Hessian of the cubic polynomial
defining the surface. We refer to the desingularization
of the corresponding nodal quartic surface as the Hessian
K3 surface of the cubic surface.

In this paper we study these Hessian K3 surfaces.
The Hessian, with its polarization given by the natural map to
$\PP^3$, determines the cubic surface. This gives a birational
isomorphism between the moduli space of $M$-lattice polarized K3 surfaces,
where $M$ is the N\'eron-Severi group of the general Hessian,
and the moduli space of cubic surfaces.

As a first step in the study of the moduli space of these K3 surfaces
we determine the transcendental lattices of various types of
Hessians. We find that certain naturally defined divisors
in the moduli space
(Hessians of cubics with a node, an Eckardt point, without a Sylvester form,
Kummer surfaces) are Heegner divisors, that is, their
transcendental lattices are sublattices of a lower rank of the
transcendental lattice of a general Hessian. Various very special cubic
surfaces (with many nodes and/or Eckardt points) have transcendental lattices
of rank two, these are usually referred to as singular
(or, nowadays, as attractive) K3 surfaces. As there is no general method
to determine the transcendental lattice of a K3 surface,
we have to use ad-hoc methods.
To study the Hessians, we use the Sylvester form of a cubic surface
which is also very useful in describing the moduli space.

In the final section of this paper we recall the classical description of the
moduli space of cubic surfaces using 19th century invariant theory.
The moduli space is a weighted projective space. We describe various
divisors and subvarieties in the moduli space. We also describe its
singular locus.

This paper is based on the PhD thesis of the first author under the direction
of the second author.

\section{The Sylvester form and the Hessian}

\subsection{The Sylvester form}\label{sylvester}
Given a general homogeneous polynomial $F$ (briefly, a form) of degree three
in four variables, there are five linear forms $x_i$ in four variables
such that any four of the five $x_i$ are linearly independent and such that
$$
\sum_{i=0}^4 x_i=0,\qquad F=\sum_{i=0}^4\lambda_ix_i^3.
$$
We refer to these equations as the Sylvester form of $F$.

The $x_i$ are uniquely determined by $F$ (up to permutation and
multiplication by a common non-zero scalar) and the $\lambda_i$ are
uniquely determined by $F$ and the $x_i$ (\cite{S}, pp.\ 125-137,
\cite{Rod}, pp.\ 72ff.\ and also \cite{C} pp.\ 295 ff.).

\subsection{The pentahedron}\label{penta}
In case $F$ has a Sylvester form, the union
of the five planes in $\PP^3$
defined by $x_i=0$ is called the pentahedron of $F$.
Each plane is called a face,
the intersection of two distinct faces is
called an edge and the intersection of three distinct faces is a
vertex of the pentahedron.
The ten edges and ten vertices are denoted by
$$
L_{ij}:\quad  x_i=x_j=0,\qquad\qquad
P_{ijk}:\quad x_i=x_j=x_k=0.
$$

\subsection{The discriminant divisor}\label{discriminant}
It is not hard to verify that the cubic surface defined by
a Sylvester form is smooth iff,
for all choices of signs, we have:
$$
\sum_{i=0}^4 \pm \sqrt{\lambda_0\ldots \hat{\lambda_i}\ldots\lambda_4}\,\neq 0.
$$
Assuming the first sign to be $+1$, this condition is equivalent to the
product of the $2^4$ such expressions, with all possible signs,
to be non-zero.
This product is homogeneous of degree 32 in in the $\lambda_i$ (in fact
the `discriminant' divisor which parametrizes singular cubic surfaces
in $\PP^{19}=\PP H^0(\PP^3,\cO(3))$ has degree 32,
cf.\ \cite{S}, App.\ III p.\ 179, \cite{GKZ}, p.\ 38).

For a general cubic form $F$ we can use the linear forms
$x_0,\ldots,x_3$ as coordinates on $\PP^3$ and $x_4=-(x_0+\ldots+x_3)$.
The general cubic surface defined by $F=0$ is now determined by
$p=(\lambda_0:\ldots:\lambda_4)\in\PP^4$, and $p$ is
unique up to permutation of the $\lambda_i$.

\subsection{The Hessian}
The Hessian $H_G$ of a cubic form $G$ in $n+1$ variables $u_i$
is the homogeneous polynomial of degree $n+1$ defined by
$$
H_G=\det\left(\frac{\partial^2 G}{\partial u_i\partial u_j}\right).
$$

\subsection{The Hessian K3 surface}\label{hesk3}
The Hessian surface of a cubic surface $S$ defined by
$F=0$ is the quartic
surface $Y$ defined by $H_F=0$ (in case $H_F\equiv 0$ we do not define the
Hessian surface).
We will call $Y=Y_S$ the Hessian of $S$.

In case $S$ is defined by a Sylvester form,
its Hessian $Y$ is defined by
$$
Y:\quad
H_F=\sum_{i=0}^4 \lambda_0\ldots \hat{\lambda_i}\ldots\lambda_4
x_0\ldots \hat{x_i}\ldots x_4=
(\lambda_0\ldots\lambda_4)(x_0\ldots x_4)\sum_{i=0}^4 \frac{1}{\lambda_ix_i}=0,
\qquad\sum_{i=0}^4 x_i=0.
$$
The ten edges $L_{ij}$ of the pentahedron
lie on the Hessian and the ten vertices $P_{ijk}$
are singular points of the Hessian.
One verifies that $S$ is smooth iff the singular locus of $Y$ consists
of these 10 points.

The desingularization $X$ of $Y$ is a K3 surface, which we refer to as the
Hessian K3 surface of $S$ (or simply Hessian if no confusion is likely).
The natural map
$$
\pi:X\longrightarrow Y,\qquad \pi(N_{ijk})=P_{ijk}
$$
contracts 10 smooth
rational curves $N_{ijk}$ to the singular points $P_{ijk}$ of $Y$ and
is an isomorphism on the complement. The strict transforms in $X$
of the lines $L_{ij}$ in $Y$ are denoted by $N_{ij}$.

\subsection{The N\'eron-Severi group of the general Hessian K3 surface}
\label{nsgen}
The ten $(-2)$-curves $N_{ijk}$ are disjoint, and so are
the ten $(-2)$-curves $N_{ij}$.
Moreover, $N_{ij}N_{klm}=1$ iff $\sharp\{i,j,k,l,m\}=3$
and is zero otherwise. A machine computation shows
that the rank of the $20\times 20$ matrix of intersection numbers of the
$N_\alpha$ is equal to $16$. A $\ZZ$-basis of the sublattice $NS_{gen}$
of $NS(X)$
generated by these curves is given by all twenty curves except
$N_{234}$, $N_{14}$, $N_{23}$ and $N_{24}$.
The discriminant of $NS_{gen}$ is equal to $2^4\cdot 3$.
From \cite{DK} it follows that $NS_{gen}$ is the
N\'eron-Severi group of a general Hessian K3 surface (i.e.\ one
with Picard number $16$). The perpendicular of $NS_{gen}$ in $H^2(X,\ZZ)$
is the transcendental lattice $T_{gen}$ which was determined in \cite{DK}
(for an alternative method, see Proposition \ref{Xks}):
$$
T_{gen}=U\oplus U(2)\oplus A_2(-2).
$$

The following lemma is very useful for determining the transcendental
lattices of `singular' Hessian K3 surfaces, that is of those with Picard
number $20$.

\subsection{Lemma}\label{slh}
Let $T$ be an even lattice of rank two:
$$
T=\left(\begin{array}{cc}
2n&a\\a&2m\end{array}\right).
$$
Then there is a primitive embedding $T\hookrightarrow T_{gen}$ if and
only if at least one among $a$, $n$ and $m$ is even. In this case
$T$ embeds into $U\oplus U(2)$.

\ts
An embedding $T\hookrightarrow T_{gen}$, is equivalent
to giving
vectors $x=(x_1,x_2,x_3,x_4,x_5,x_6)$,
$y=(y_1,y_2,y_3,y_4,y_5,y_6)\in \ZZ^6=T_{gen}$, such that
$$
\begin{array}{ll}
2x_1x_2+4x_3x_4-4(x_5^2+x_6^2-x_5x_6)&=2n,\\
2y_1y_2+4y_3y_4-4(y_5^2+y_6^2-y_5y_6)&=2m,\\
x_1y_2+x_2y_1+2(x_3y_4+x_4y_3)-2(2x_5y_5+2x_6y_6-x_5y_6-x_6y_5)&=a.
\end{array}
$$

If $n\equiv m\equiv a\equiv 1(2)$, this system has no solution,
since from the first two equations we get $x_1x_2\equiv 1(2)$ and
$y_1y_2\equiv 1(2)$, so $x_1\equiv x_2\equiv y_1\equiv y_2\equiv
1(2)$ and hence $x_1y_2+x_2y_1\equiv 0(2)$. Thus in the third equation
the first member is even, while the second is odd.

Conversely, let $x:=(n,1,0,0,0,0)$, then $y$ must satisfy
$$
\begin{array}{ll}
y_1y_2+2y_3y_4-2(y_5^2+y_6^2-y_5y_6)&=m,\\
ny_2+y_1&=a.
\end{array}
$$
Substituting $y_1=a-ny_2$ in the first equation, we get:
$$
2y_3y_4-2(y_5^2+y_6^2-y_5y_6)=ny_2^2-ay_2+m.
$$
Hence we can choose:
$$
y=(a-ny_2,y_2,1,(ny_2^2-ay_2+m)/2,0,0),
$$
unless $n\equiv a(2)$ and $m\equiv 1(2)$,
the only case in which $ny_2^2-ay_2+m$ is odd for all $y_2\in\ZZ$.
Since, by assumption, one of $n$, $m$ and $a$ is even, in this
case $n$ (and $a$) must be even. Now we begin with
$y=(1,m,0,0,0,0)$ and
the same argument provides an $x$.
It is easy to see that with these choices of $x$ and $y$ the embedding
of $T$ is primitive.
\qed

\section{Eckardt points}

\subsection{Definition} A smooth cubic surface contains 27 lines and has
45 plane sections which are unions of three lines.
In case three coplanar lines
meet in a single point, that point is called an Eckardt point.

\subsection{The Sylvester form and Eckardt points}
A smooth cubic surface $S$ defined by a Sylvester form
$\sum \lambda_ix_i^3$ as in
\ref{sylvester} has an Eckardt point iff $\lambda_i=\lambda_j$ for some
$i\neq j$ (\cite{S}, pp.\ 145 ff.)

If $\lambda_i=\lambda_j$, the corresponding Eckardt point $P$ of $S$
is the vertex $P_{klm}$ of the
pentahedron (with $\{i,\ldots,m\}=\{0,\ldots,4\}$). The plane $x_i+x_j=0$
(which is the tangent plane $T_PS$ to $S$ in $P_{klm}$) cuts out the three
lines meeting in $P_{klm}$
and $S$ has an involution (induced by permuting $x_i$ and $x_j$).
In particular, $S$ has an Eckardt point iff a vertex of the
pentahedron lies on $S$.

A Sylvester form with $\lambda_0=\lambda_1=\lambda_2$ defines
a cubic surface with three Eckardt points in general, there are
$4$ Eckardt points if also $\lambda_3=\lambda_4$ and there are 10
Eckardt points if all $\lambda_i$ are equal (this surface is known as
the Clebsch diagonal surface).
In case $\lambda_0=\lambda_1$ and $\lambda_2=\lambda_3$, but otherwise the
$\lambda_i$'s are distinct, the surface has $2$ Eckardt points.
In case $\lambda_0=\lambda_1=\lambda_2=\lambda_3$  but
$\lambda_0\neq \lambda_4$, the surface has six Eckardt points.

\subsection{New curves on the Hessian}
For an Eckardt point $P$
the intersection of the tangent plane $T_PS$ and the Hessian $Y_S$
consists of the line $L_{ij}$ (with multiplicity two) and
of the pair of lines defined by
$$
\lambda_k\lambda_lx_kx_l+\lambda_k\lambda_mx_kx_m+\lambda_l\lambda_mx_lx_m=0
,\qquad
x_k+x_l+x_m=0,
$$
which meet in $P_{klm}$, we refer to these as new lines.
The strict transform of these lines in the Hessian K3 surface are two new
disjoint $(-2)$-curves which both meet $N_{klm}$, but none of the other
$N_{abc}$'s.

\subsection{Clebsch diagonal cubic surface}\label{x10}
The diagonal cubic surface is the smooth cubic surface $S_{10}$ defined by
the Sylvester form:
$$
S_{10}:\qquad\sum_{i=0}^4 x_i^3=0,\qquad \sum_{i=0}^4 x_i=0,
$$
so all $\lambda_i=1$. This surface has $10$ Eckardt points and thus
its Hessian $Y_{10}$ has $10$ new pairs of lines.
Let $\omega$ be a primitive cube root of unity and let
$$
C_{234}=\{(s:-s:t:t\omega:t\omega^2):\;(s:t)\in\PP^1\}\qquad
\subset Y_{10}\qquad (\omega^2+\omega+1=0)
$$
be one of the new lines on the Eckardt point
$P_{234}=(1:-1:0:0:0)$. For $\sigma\in S_5$, the symmetric group on
$\{0,\ldots,4\}$, we denote by $C_{\sigma(2),\sigma(3),\sigma(4)}$
the line obtained by permuting the coordinates. In particular,
$$
C_{234}=C_{342}=C_{423},\qquad C_{243}=C_{324}=C_{432}=
\{(s:-s:t:t\omega^2:t\omega):\;(s:t)\in\PP^1\}
$$
is the other new line on $P_{234}$. We use the same names to denote
the corresponding $(-2)$-curves on the Hessian K3 surface $X_{10}$ of
$S_{10}$.

The intersection numbers of $N_{ab}$ and $N_{abc}$ with $C_{klm}$
are zero unless $\sharp\{a,b,k,l,m\}=5$ resp.\
$\sharp\{a,b,c,k,l,m\}=3$, in which case the intersection number is $1$.
Moreover,
$C_{abc}$ and $C_{klm}$ intersect, with intersection number 1,
iff there are two equal indices in the same order
(up to cyclic permutations,
e.g.\ $C_{124}$ intersects $C_{041}$), $C_{abc}$ and $C_{abd}$
meet in a point on $X_{10}$ whose image in
$Y_{10}$ has coordinates $x_a=\omega$, $x_b=\omega^2$,
$x_c=x_d=1$ and $x_e=-1$ where $\{a,b,c,d,e\}=\{0,1,2,3,4\}$.
The following classes in the N\'eron-Severi group $NS_{10}$
of $X_{10}$
$$
c_{ij}=C_{abc}-C_{acb}\quad\in NS_{10}:=NS(X_{10}),\qquad
\{i,j,a,b,c\}=\{0,\ldots,4\},
$$
are easily seen to be perpendicular to $NS_{gen}$. They satisfy:
$$
c_{ij}^2=-4,\qquad c_{ik}\cdot c_{il}=2,\qquad c_{ij}\cdot c_{kl}=0
$$
for distinct indices $i,j,k,l$, so we have a copy of $A_4(-2)$ in $T_{gen}$
(cf.\ Remark \ref{a4}).
The symmetric group $S_5$ is a subgroup of $Aut(X_{10})$
(acting via permutations of the coordinates) and acts via permutation
of the indices on the $c_{ij}$, so it acts as $W(A_4)\cong S_5$
on the copy of $A_4(-2)$.

\subsection{Proposition}\label{Xks}
Let $X_k$ be the Hessian K3 surface of a smooth
cubic surface with Sylvester form having
$k$ Eckardt points (in particular, $k\in\{0,1,2,3,4,6,10\}$).
Let $T_k$ be the transcendental lattice of $X_k$.
$$
\begin{array}
{ll}
T_0=T_{gen}=U\oplus U(2)\oplus A_2(-2)\qquad &{\rm discr}(T_{gen})=48,\\
T_1=U\oplus U(2)\oplus <-12>\qquad &{\rm discr}(T_1)=-48,\\
T_2=U\oplus <4>\oplus <-12>\qquad & {\rm discr}(T_2)=48,\\
T_3=U\oplus U(6)\qquad &{\rm discr}(T_3)=36,\\
T_4=U(3)\oplus <4> \qquad& {\rm discr}(T_4)=-36,\\
T_6=U\oplus <24> \qquad& {\rm discr}(T_6)=-24.\\
T_{10}=\left(\ZZ^2,\left(\begin{array}{cc}4&1\\1&4\end{array}\right)\right)
\qquad& {\rm discr}(T_{10})=15.
\end{array}
$$

\ts
We start with the computation of $NS_{10}$ and $T_{10}$.
The $40\times 40$ matrix of intersection products of the curves
$N_{\alpha}$ and $C_\beta$ has rank 20, hence the N\'eron-Severi group
of $X_{10}$ has maximal rank. A $\ZZ$-basis of the lattice generated by
these 40 curves is given by the basis of $NS_{gen}$ (given in \ref{nsgen})
and the 4 curves $C_{234}$, $C_{134}$, $C_{124}$ and $C_{032}$.
A machine verified that the discriminant of this lattice is $-15$.
Since $15$ is square free, we can conclude that $NS_{10}$ is generated
by these 40 curves. It is easy to find two orthogonal $E_8(-1)$'s in this
lattice, for example:
$$
\begin{array}{ccccccccccccc}
N_{034}&--&N_{04}&--&N_{024}&--&C_{024}&--&C_{124}&--&
C_{041}&--&N_{23}\\
&&&&\mid&&&&&&&&\\
&&&&N_{24}&&&&&&&&
\end{array}
$$
$$
\begin{array}{ccccccccccccc}
N_{12}&--&N_{012}&--&N_{01}&--&C_{234}&--&C_{142}&--&
C_{143}&--&N_{134}\\
&&&&\mid&&&&&&&&\\
&&&&N_{013}&&&&&&&&
\end{array}
$$
Using a machine, it is not hard to find the rank 4 lattice $L'$
such that $NS_{10}\cong E_8(-1)^2\oplus L'$,
we found a copy of $U$ in $L'$ (a standard basis is:
$-2N_{134}-5N_{124}-N_{123}+4N_{034}+5N_{024}+N_{023}+7N_{014}+6N_{013}
+8N_{012}+13N_{01}+2N_{02}+8N_{04}+N_{12}-N_{13}+5C_{234}-3C_{134}
-4C_{124}+C_{032}$ and
$-2N_{134}-5N_{124}-N_{123}+4N_{034}+5N_{024}+N_{023}+7N_{014}+6N_{013}
+8N_{012}+13N_{01}+2N_{02}+8N_{04}+N_{12}-N_{13}+5C_{234}-3C_{134}
+4C_{124}$).
The orthogonal complement of $U$ in $L'$ is isomorphic to the lattice
denoted by $T_{10}(-1)$ in the theorem
(a basis for it is:
$2N_{124}+2N_{123}-5N_{034}-2N_{024}-2N_{023}-N_{014}+2N_{012}+N_{01}-3N_{03}
-4N_{04}+3N_{12}+2N_{13}-3N_{34}+C_{234}+C_{134}-C_{032}$ and
($-4N_{134}-9N_{124}-4N_{123}+8N_{034}+9N_{024}+4N_{023}+11N_{014}
+9N_{013}+11N_{012}+19N_{01}+4N_{02}+2N_{03}+14N_{04}-N_{12}-3N_{13}
+7C_{234}-5C_{134}-6C_{124}+2C_{032}$).
In general, the embedding of a rank 20
lattice like $NS_{10}$ into the K3 lattice $E_8(-1)^2\oplus U^3$
is not unique up to isomorphism. However, in this case we know that
the orthogonal complement of $NS_{10}$ in $H^2(X_{10},\ZZ)$ is an
even rank two lattice with discriminant 15. The classification
of positive definite binary quadratic forms shows that there
are two such lattices, one is $T_{10}$,
the other has matrix with rows $2,1$ and $1,8$.
The discriminant form of the N\'eron-Severi lattice
is the opposite of the one of the transcendental lattice
(after identifying the discriminant groups).
The discriminant form $NS_{10}$ is equal to minus the one of $T_{gen}$,
but the discriminant form of the other lattice is different,
hence we found $T_{10}$.

A result of Nikulin implies that the embedding of $T_{10}$ in
the K3 lattice is unique up to isometry (cf.\ \cite{M}, cor.\ 2.10).
Choosing a convenient embedding, it is then easy to find an explicit
isomorphism of $NS_{10}$ with a sublattice of the K3 lattice.

To find the other lattices, we use that the number of moduli $d_k$ of K3
surfaces with a fixed number $k$ of Eckardt points is easy to find using the
Sylvester form. By the Torelli theorem
the rank of the transcendental lattice
of a general K3 in this family is at least $d_k+2$, hence
$NS_k$ has rank at most $20-d_k$.
Computations show that in each case the lattice spanned
by $NS_{gen}$ and the $C_\alpha$'s on the general member of the family
has rank $20-d_k$.
It is not hard to check that this lattice is perpendicular
to certain elements in $NS_{10}$,
choosing the Sylvester forms defining the $X_k$ such that
$\lambda_0=\lambda_1=\lambda_2=\lambda_3$ for $X_6$,
$\lambda_0=\lambda_1=\lambda_2$ and $\lambda_3=\lambda_4$ for $X_4$,
$\lambda_0=\lambda_1=\lambda_2$ for $X_3$, $\lambda_0=\lambda_1$ and
$\lambda_2=\lambda_3$ for $X_2$ and $\lambda_0=\lambda_1$ for $X_1$
we have:

\begin{enumerate}
\item $NS_6=<c_{04}-c_{14}+c_{24}-c_{34}>^{\bot}\hookrightarrow NS_{10}.$

\item
$NS_4=<c_{14}-c_{24}+c_{03}-c_{13}+c_{23} >^{\bot} \hookrightarrow NS_{10};$

\item $NS_3=<c_{34}>^{\bot}\hookrightarrow NS_4$ and
\noindent
$NS_3=<c_{23}-c_{13}+c_{03}>^{\bot}
\hookrightarrow NS_6;$

\noindent Obviously, to embed $NS_2$ in $NS_6$ we need to make a
permutation of the indices and consider $X'_2$, i.e.\ the Hessian
of the cubic with two Eckardt points and Sylvester representation
with coefficients $\lambda_0=\lambda_1$ and $\lambda_2=\lambda_3$.
Thus the embedding is:
$NS'_2=<c_{02}-c_{03}-c_{12}+c_{13}>^{\bot}\hookrightarrow NS_6;$

\item $NS_2=<c_{02}-c_{12}>^{\bot}\hookrightarrow
NS_4$;

\item $NS_1=<c_{34}>^{\bot}\hookrightarrow NS_2$ and
\noindent
$NS_1=<c_{02}-c_{12}>^{\bot}\hookrightarrow
NS_3;$

\item $NS_{Hess}=<c_{01}>^{\bot}\hookrightarrow NS_1$.
\end{enumerate}
Next we determined the transcendental lattice $T_i$
as $NS_i^\perp$ in the K3 lattice.
In particular, in this way we determined $NS_{gen}$ and its
perpendicular $T_{gen}$ in the K3 lattice.
\qed

\subsection{Remark}\label{a4}
We give an explicit primitive embedding of $A_4(-2)$
in $T_{gen}$. Recall that $T_{gen}=U\oplus U(2)\oplus A_2(-2)$,
an element $x\in T_{gen}$ will be written as $x=(x_1,\ldots,x_6)\in\ZZ^6$,
the quadratic form on $T_{gen}$ is $2x_1x_2+4x_3x_4-2(2x_5^2-x_5x_6+2x_6^2)$.

A basis of $A_4(-2)$ is given by the $c_{ij}$ (cf.\ section \ref{x10}) and
the embedding of $A_4(-2)$ in $T_{gen}$ is:
$$
c_{01}= (0, 0, 0, 0, 0, 1),\qquad
c_{12}= (0, 0, 0, -2, 1, 0)
$$
$$
c_{23}=(0, 0, 1, 2, -2, -1),\qquad
c_{34}=(4, -2, -3, -5, 4, 2).
$$
The orthogonal complement of this sublattice in $T_{gen}$ is isomorphic
to $T_{10}$ and is spanned by
$$
t_1=(2, 1, 0, 0, 0, 0),\qquad t_2=(5, -2, -3, -6, 4, 2),\qquad
{\rm with}\quad
(t_i\cdot t_j)=\left(\begin{array}{cc}4&1\\1&4\end{array}\right).
$$
Moreover, the lattice
$$
T_{10}\oplus A_4(-2)\;\hookrightarrow\; T_{gen}=U\oplus U(2)\oplus A_2(-2)
$$
has index 5 and $T_{gen}$ is generated by $T_{10}\oplus A_4(-2)$
and
$$
\mbox{$\frac{1}{5}$}
(2(t_1+t_2)+c_{01}+2c_{12}+3c_{23}+4c_{34}))=(6, -2, -3, -6, 4, 2).
$$

\section{Nodal cubic surfaces}

\subsection{Cayley's four nodal cubic surface}
A cubic surface can have up to four nodes. There is a unique surface $S_{4n}$
having four nodes, it is Cayley's cubic surface defined by
$$
\frac{1}{z_0}+\frac{1}{z_1}+\frac{1}{z_2}+\frac{1}{z_3}=0,
$$
it has the Sylvester form:
$$
S_{4n}:\qquad\sum x_i=0,\qquad x_0^3+x_1^3+x_2^3+x_3^3+\frac{1}{4}x_4^3=0.
$$
Its nodes (in the Sylvester form) are $p_0=(-1:1:1:1:-2)$ and
the other three are obtained by permuting the first four coordinates.

These points are also singular on the Hessian $Y_{4n}$, and they give rise
to $(-2)$-curves $M_i$ on the Hessian K3 surface $X_{4n}$.
The lines $L_{ij}$ connecting two nodes $p_i,p_j$
of the Hessian are on the Hessian. Thus we find another six $(-2)$-curves
on $X_{4n}$ which we also denote by $L_{ij}$.

In case the cubic surface has $k$ nodes, we denote the $(-2)$-curves on
the Hessian K3 by $M_1,\ldots,M_k$ and $L_{ij}$ ($1\leq i<j\leq k$).

\subsection{Proposition}\label{Xkns}
Let $T_{kn}$ be the transcendental lattice of the Hessian K3 surface
of a general cubic surface with $k$ nodes. Then we have:
$$
T_{4n}=<2>\oplus <6>
$$
and:
$$
T_{3n}=T_{4n}\oplus <-2>,\qquad
T_{2n}=T_{4n}\oplus <-2>^2,\qquad
T_{1n}=T_{4n}\oplus<-2>^3.
$$

\ts
The proof is similar to the one of Proposition \ref{Xks}. The
non-zero intersection
numbers involving the ten new $(-2)$-curves on the Hessian K3 surface
$X_{4n}$ of
the Cayley cubic are:
$$
\begin{array}{lll}
(M_k,L_{ij})=1\quad & \textrm{if} \quad & \#\{i,j,k\}=2,\\
(N_{ij4},L_{kl})=1\quad &\textrm{if}& \#\{i,j,k,l\}=4,\\
(N_{ij},L_{kl})=1 & \textrm{if}& \#\{i,j,k,l\}=2,
\end{array}
$$
and of course $M_i^2=L_{ij}^2=-2$. A machine computation confirms that
the rank of the lattice spanned by the 30 $(-2)$-curves in $X_{4n}$
is indeed 20 and that a basis of this lattice is given by the
16 curves in the basis of $NS(X_{gen})$ (cf.\ \ref{nsgen}) and the 4 curves
$L_{01},L_{02},L_{03},L_{12}$. The discriminant of this lattice is $-12$.
Thus the discriminant of the N\'eron-Severi lattice is either $-12$ or
$-3$, but if it were $-3$ the transcendental lattice had to be $A_2$,
which is impossible by Lemma \ref{slh}. Thus these 20 curves span the
N\'eron-Severi lattice.

There are two perpendicular $E_8(-1)$'s in the N\'eron-Severi lattice:
$$
\begin{array}{ccccccccccccc}
N_{012}&--&N_{02}&--&N_{024}&--&N_{24}&--&N_{234}&--&
N_{34}&--&N_{134}\\
&&&&\mid&&&&&&&&\\
&&&&N_{04}&&&&&&&&
\end{array}
$$
$$
\begin{array}{ccccccccccccc}
N_{013}&--&N_{03}&--&L_{03}&--&M_{3}&--&L_{23}&--&
M_{2}&--&L_{12}\\
&&&&\mid&&&&&&&&\\
&&&&M_{0}&&&&&&&&
\end{array}
$$
A further computation then shows that
$$
NS(X_{4n})\cong E_8(-1)^2\oplus U\oplus <-2>\oplus <-6>.
$$
Thus $T_{4n}$ has rank two, discriminant $12$ and its
discriminant form is the opposite of the one on $NS(X_{4n})$ and there
is a unique such lattice, which is $<2>\oplus <6>$.

Another way to obtain this result is by realizing the Hessian as a
double cover of the plane by projecting from the node $M_0$.
This double cover branches over the union of three lines and a cubic
with a node.
The pencil of lines on the node of the cubic curve defines
an elliptic fibration on $X_{4n}$ which turns out to have
six singular fibers, three of type $I_6$ and three of type $I_2$.
According to the table in \cite{ShZh} there is only one such K3 surface
(case 4) and its transcendental lattice is indeed $<2>\oplus <6>$.
This surface also appears as a double cover of the plane branched over
three conics in \cite{Persson}, Example 3, p.304.

There is a unique, up to isomorphism, embedding of $T_{4n}$
into the K3 lattice. It is not hard to find an
explicit isomorphism of the N\'eron-Severi lattice with the orthogonal
complement of the image of $T_{4n}$.
The N\'eron Severi lattice $NS_{kn}$ of $X_{kn}$ is now
easily seen to be:
$$
NS_{3n}=(T_{4n}
\oplus <M_3>)^{\bot},\quad\ldots\quad,
NS_{1n}=(T_{4n} \oplus
<M_1,M_2,M_3>)^{\bot},
$$
next one finds $T_{kn}$ as the orthogonal complement of $NS_{kn}$.
\qed

\section{Cubics with Eckardt points and nodes}

\subsection{}
In this section we consider the Hessians $X_{anb}$
of three cubic surfaces with $a$ nodes and $b$ Eckardt points
with $(a,b)=(1,6),(1,4)$ and $(3,4)$. These K3 surfaces
are `singular', i.e.\
have N\'eron-Severi groups of rank 20. In two of the three
cases we succeeded in determining the transcendental lattice completely.

\subsection{Proposition}
\begin{enumerate}
\item{}
The Hessian K3 surface $X_{1n6}$
of the cubic surface $S_{1n6}$ with one node and 6 Eckardt points
with Sylvester form
$$
\sum x_i=0,\qquad x_0^3+x_1^3+x_2^3+x_3^3+\frac{1}{16}x_4^3=0
$$
has transcendental lattice
$$
T_{1n6}=\,<2>\oplus <24>.
$$
\item{}
The Hessian K3 surface $X_{1n4}$
of the cubic surface $S_{1n4}$ with one node and 4 Eckardt points
with Sylvester form
$$
\sum x_i=0,\qquad x_0^3+x_1^3+x_2^3+\frac{4}{9}x_3^3+\frac{4}{9}x_4^3=0
$$
has transcendental lattice which is either
$$
<6>\oplus <12>\qquad{\rm or}\qquad <2>\oplus <4>.
$$
\item{}
The Hessian K3 surface $X_{3n4}$
of the cubic surface $S_{3n4}$ with three nodes and 4 Eckardt points
with Sylvester form
$$
\sum x_i=0,\qquad x_0^3+x_1^3+x_2^3+4x_3^3+4x_4^3=0
$$
has transcendental lattice
$$
T_{3n4}=\,<4>\oplus <6>.
$$
\end{enumerate}

\ts
The only singular point of $S_{1n6}$ is the point $p_0=(1:\ldots:1:-4)$
which is also singular on the Hessian, which has eleven singular points.
This point does not lie on any of the 10 lines of the Hessian,
so the N\'eron-Severi group of $X_{1n6}$ contains the
direct sum of the N\'eron-Severi group of the general Hessian K3
with 6 Eckardt points and the lattice $<-2>$, spanned by the $(-2)$-
curve over $p_0$:
$$
NS_6\oplus <-2>\,=\,E_8(-1)^2\oplus U\oplus <-24>\oplus<-2>
\;\hookrightarrow\;NS_{1n6}.
$$
As the discriminant of this sublattice is $-48$, the discriminant
of $NS_{1n6}$ is either $-48$ or there is a vector $v$, primitive
in $NS_6\oplus <-2>$, such that $v/2\in NS_{1n6}$.
As $E_8(-1)^2\oplus U$ is a direct summand, we can write $v=ae+bf$
with $a,b\in\{0,1\}$ and $e^2=-24,ef=0,f^2=-2$. As
$(v/2)^2=-6a^2+b^2/2\in2\ZZ$ (since $NS_{1n6}$ is even) we get $b=0$,
as $<-24>$ is primitive in $NS_{6}$ and hence, after specialization, in
$NS_{1n6}$, we get $a=0$. Therefore $NS_6\oplus <-2>\cong NS_{1n6}$.
The transcendental lattice $T_{1n6}$ is a positive definite
even rank two lattice whose discriminant form is the opposite of the one
of $NS_{1n6}$, and there is only one such lattice.

The unique singular point of $S_{1n4}$ is $p_0=(-2:-2:-2:3:3)$
which does not lie on any of the lines in the Hessian.
Thus one finds as above an embedding of lattices:
$$
NS_4\oplus <-2>\,=\,E_8(-1)^2\oplus U(3)\oplus <-4>\oplus<-2>
\;\hookrightarrow\;NS_{1n4}.
$$
In this case we cannot decide whether the sublattice has index 1 or 3.
In case the index is one, one has to use the discriminant form to find
$T_{1n4}$, if the index is three then the discriminant of
$T_{1n4}$ is $8$, and there is a unique such even positive definite lattice.

The three singular points of $S_{3n4}$ are:
$$
p_0=(-2:2:2:-1:-1),\qquad
p_1=(2:-2:2:-1:-1),\qquad
p_2=(2:2:-2:-1:-1).
$$
The lines $L_{ij}$ spanned by $p_i$ and $p_j$ lie on the Hessian.
The $p_i$'s do not lie on any of the 10 lines $N_{ij}$ but
the lines $L_{ij}$ and $N_{ij}$ intersect, finally $L_{ij}$ and $N_{klm}$
intersect if $\sharp \{i,j,k,l,m\}=5$.
Of the lines $C_{abc}$ on the Hessian of a cubic with Eckardt points,
only the lines $C_{012}$ and $C_{021}$ remain, the others collapse
onto the $L_{ij}$.
A machine computation shows that the lattice spanned by these
$(-2)$-curves has rank $20$ and discriminant $-24$.
A basis is given by the 16 curves in the basis of $NS(X_{gen})$
(cf.\ \ref{nsgen}) and $C_{012},L_{01},L_{02},L_{12}$.
Since the determinant of an even rank two lattice cannot be $6$,
we conclude that this lattice is $NS_{gen}$. To find its discriminant form
we used the following two copies of $E_8(-1)$:
$$
\begin{array}{ccccccccccccc}
N_{23}&--&N_{123}&--&N_{13}&--&N_{013}&--&N_{01}&--&
L_{01}&--&M_{0}\\
&&&&\mid&&&&&&&&\\
&&&&N_{134}&&&&&&&&
\end{array}
$$
$$
\begin{array}{ccccccccccccc}
N_{124}&--&N_{24}&--&N_{024}&--&N_{04}&--&N_{034}&--&
L_{12}&--&M_2\\
&&&&\mid&&&&&&&&\\
&&&&N_{02}&&&&&&&&
\end{array}
$$
and we found an isomorphism
$$
NS_{3n4}\cong E_8(-1)^2\oplus U\oplus <-4>\oplus <-6>.
$$
There is a unique even, rank 2, positive definite lattice with
opposite discriminant group.
\qed

\section{Cubic forms without a Sylvester form}\label{nosylvester}

\subsection{No Sylvester form} We defined the Sylvester form
of a cubic form $F$
in 4 variables using 5 linear forms $x_0,\ldots,x_4$ in 4 variables
such that any four of the five are linearly independent, cf.\ \ref{sylvester},
normalized such that $\sum x_i=0$ and satisfying $F=\sum \lambda_ix_i^3$
for some $\lambda_i\in\CC$.

There are basically two types of cubic forms defining smooth cubic surfaces
which do not admit a Sylvester form (we do not know precise results on
forms defining singular surfaces). In the first case one still has
$F=\sum \lambda_ix_i^3$,
but at least four of the five $x_i$ are linearly
dependent, say $x_4=a_1x_1+a_2x_2+a_3x_3$. Then one has
$F=\lambda_0x_0^3+G(x_1,x_2,x_3)$ with a homogeneous polynomial $G$
in three variables.
The corresponding cubic surfaces are called cyclic surfaces. Note that the
Hessian of a cyclic cubic surface is reducible, which proves that these
surfaces do not have a Sylvester form.

The second type of cubic forms $F$ which do not admit a Sylvester forms
are those which are not a sum of five cubes of linear forms. These are
obtained as a limit of Sylvester forms where at least two
of the five $x_i$ coincide, we will refer to these as non-Sylvester forms.
The corresponding surfaces were described in \cite{S} and \cite{Rod}.
They are of two types, denoted here by $ns1$ and $ns2$,
the dimensions of these families are three and two respectively.
We give explicit one parameter families of Sylvester forms
specializing to these surfaces in the proof of Theorem \ref{nsinv}.

\subsection{The first case}\label{ns1}
A cubic surface of type $ns1$ is  defined by a form $F$ of the following type:
$$
S_{ns1}:\qquad F=x_1^3+x_2^3+x_3^3-x_0^2(a_0x_0+3a_1x_1+3a_2x_2+3a_3x_3)=0.
$$
This form defines a smooth cubic surface iff, for all choices of sign,
$$
a_0+2(\pm a_1^{3/2}\pm a_2^{3/2}\pm a_3^{3/2})\neq 0.
$$
In case one of the $a_i=0$ for $1\leq i\leq 3$, the form $F$ defines a cyclic
surface. We will assume that $S_{ns1}$ is smooth and that
$a_1a_2a_3\neq 0$ from now on.

The Hessian of $F$ is, up to a scalar multiple,
$$
H_F=x_1x_2x_3(a_0x_0+\ldots+a_3x_3)+x_0^2(a_1^2x_2x_3+a_2^2x_1x_3+a_3^2x_1x_2).
$$
Obviously, $p_0=(1:0:0:0)$ is singular on the Hessian surface $Y$.
Projecting from $p_0$ realizes the Hessian surface as a double cover of the
plane (with coordinates $x_1,x_2,x_3$) branched over the curve defined by:
$$
x_1x_2x_3[a_0^2x_1x_2x_3-
4(a_1x_1+a_2x_2+a_3x_3)(a_1^2x_2x_3+a_2^2x_1x_3+a_3^2x_1x_2)]=0,
$$
that is, the union of three lines and a nonsingular cubic curve.
As a singular point of the discriminant curve is the image of a singular point
or the image of a line on the Hessian surface it is now easy to determine
the singular locus of $Y$. One finds that $Y$ has exactly 7 singular points,
$p_0$ and the six points:
$$
p_1=(0:1:0:0),\quad p_2=(0:0:1:0),\quad p_3=(0:0:0:1),
$$
(these points lie on lines $x_i=x_j=0$ ($1\leq i<j\leq 3$) in $Y$ which project
to points),
$$
q_1=(0:0:a_3:-a_2),\quad q_2=(0:-a_3:0:a_1),\quad q_3=(0:a_2:-a_1:0).
$$
In particular, the cubic surface $S_0$ defined by $F$ cannot have a Sylvester
form (in that case the Hessian has exactly 10 singular points) nor is it
cyclic (in that case the Hessian surface is reducible).

\subsection{The second case}\label{ns2}
A cubic surface of type $ns2$ is defined as follows.
$$
S_{ns2}:\qquad G=x_1^3+x_2^3+2\lambda x_3^3-3x_3(\mu x_1x_3+x_2x_3+x_0^2)=0.
$$
If $\mu=0$ this is a cyclic cubic surface and the surface $S_{ns2}$ is smooth iff
for all choices of sign:
$$
\lambda\pm \mu^{3/2}\pm 1\neq 0.
$$
The Hessian of $G$ is (up to scalar multiple):
$$
H_G=x_1x_2x_3(-2\lambda x_3+\mu x_1+ x_2)+x_3^3(x_1+\mu^2x_2)-x_0^2x_1x_2.
$$
Again $p_0=(1:0:0:0)$ is a singular point. Projecting from $p_0$ gives the
branch curve defined by:
$$
x_1x_2x_3[x_1x_2(-2\lambda x_3+\mu x_1+ x_2)+x_3^2(x_1+\mu^2x_2)]=0,
$$
which is the union of three lines and a smooth cubic curve.
The Hessian surface has only four singular points, $p_0$, and
$$
p_1=(0:1:0:0),\quad p_2=(0:0:1:0),\quad q=(0:1:-\mu:0).
$$
As before, we conclude that these surfaces do not admit a Sylvester form nor
are they cyclic.

\subsection{Proposition}\label{tnsk}
Let $X_{nsk}$, $k=1,2$, be the Hessian K3 surface of a general cubic surface
$S_{nsk}$ which does not admit a Sylvester form. Let $T_{nsk}$ be the
transcendental lattice of $X_{nsk}$. Then we have:
$$
T_{ns1}= U\oplus U(2)\oplus <-4>,\qquad T_{ns2}=U\oplus U(2).
$$

\ts
First we determine the N\'eron-Severi group of $X_{ns2}$ using
the elliptic fibration defined by the pencil of lines on the image $\bar{q}$
of $q$ in $\PP^2$. From a study of the branch curve one finds that
there are 6 bad fibers, of type $I_0^*,I_8^*,I_1,I_1,I_1,I_1$.
The lines $x_i=x_3=0$ ($i=1,2$) define sections of the fibration.
The rank of the N\'eron-Severi group of the general $X_{ns2}$ is $18$
(the surfaces $X_{ns2}$ have two moduli),
which is two plus the number of components of
the bad fibers not meeting one of the sections. Hence the Mordell-Weil group
of the fibration is finite. The torsion subgroup of a fiber of type $I_0^*$ is
$(\ZZ/2\ZZ)^2$
and the two-torsion of an $I_1$-fiber is $\ZZ/2\ZZ$. As the torsion
of the Mordell-Weil group injects in the torsion subgroup of the fibers it
must be either trivial or be isomorphic to $\ZZ/2\ZZ$. As there are
two sections, we conclude that the Mordell-Weil
group of the elliptic fibration is $\ZZ/2\ZZ$.
The Shioda-Tate formula (\cite{Shiodainose}, Lemma 1.3)
now shows that the discriminant of $NS(X_{ns2})$ is
$4$.

The $13$ $(-2)$-curves in the $I_8^*$ fiber (which map to the line $x_3=0$
in $\PP^2$, note the configuration of curves is two copies of $D_6$,
mapping to $x_1=x_3=0$ and $x_2=x_3=0$, linked with a vertex which is a
$(-2)$-curve mapping to the line $x_3=0$),
the $(-2)$-curve over $\bar{q}$, the two sections and the four of the five
$(-2)$-curves in the $I_0^*$ fiber (these 4
curves map to the point $x_1=x_2=0$ in $\PP^2$) give a diagram with
20 vertices. It is a square with 3 points on each edge and each of the four
vertices of the square is
linked to a further point, which is not linked to anything else.
In particular, it is very easy to find two perpendicular $E_8$'s in the
lattice generated by these 20 curves
and to compute that the orthogonal complement to this sublattice
is $U(2)$. Hence $NS(X_{ns2})\cong E_8(-1)^2\oplus U(2)$. By \cite{Nikulin},
Thm.\ 1.14.4, there is, up to isometry, a unique embedding of $U(2)$
into $U^3$. From this one finds that
$T_{ns2}=NS(X_{ns2})^\perp\cong U(2)\oplus U$.

Similar to the case of $X_{ns2}$,
the family of lines on the image of $q_3$ in $\PP^2$ defines an elliptic
fibration $X_{ns1}$ with two sections, $x_i=x_3=0$, $i=1,2$.
There are 7 singular fibers of type
$I_4^*,I_4,I_0^*,I_1,I_1,I_1,I_1$. There are 15 components of bad fibers
which do not meet one of the sections,
so the Mordell-Weil group of this fibration
finite. As before we conclude that the Mordell-Weil
group of the elliptic fibration is $\ZZ/2\ZZ$.
The Shioda-Tate formula shows that the discriminant of $NS(X_{ns1})$ is
$2^4$.

It is amusing to observe that the
following 20 $(-2)$-curves have an intersection diagram
which is a cube with a point on the middle of each edge:
the $9$ components of the
$I_4^*$ fiber;
the following $3$ components of the $I_4$-fiber:
the line $l:\,x_0=a_1x_1+a_2x_2+a_3x_3=0$ in $Y_{ns1}$
and the $(-2)$-curves over $q_1$ and $q_2$ in $X_{ns2}$;
the $(-2)$-curves over $p_0,q_3$;
the two sections; and 4 of the 5 components of the $I_0^*$
fiber.

It is clear that $X_{ns2}$ is a limit of $X_{ns1}$ (two of the lines in the
branch curve become tangent to the cubic component, equivalently, the fibers
of type $I_4$ and $I_4^*$ coalesce to a fiber of type $I_8^*$).
Thus we get an
inclusion $NS(X_{ns1})\subset NS(X_{ns2})$ and it is not hard to see that
$NS(X_{ns1})=(n_1-n_2)^\perp$ where $n_1,n_2$ are two disjoint $(-2)$-curves
in $X_{ns2}$, they are the inverse image of the strict transform of the
exceptional divisor in the first blow up of $\PP^2$ in the image of
$p_1$ and $p_2$. Thus the lattice $T(X_{ns2})\oplus <-4>$ is a sublattice of
$T(X_{ns1})$, but since these lattices have the same discriminant, they are
equal.
\qed

\subsection{Remark}
Specializing a general Hessian to one of type $ns1$ induces an
inclusion $NS_{gen}\hookrightarrow NS_{ns1}$.
One can verify that $NS_{gen}^\perp$ in $NS_{ns1}$ is generated
by a class with selfintersection $-12$. Note that $t=(1,2)\in A_2(-2)$
has $t^2=-12$ and that its orthogonal complement is $<(1,0)>$ with
$(1,0)^2=-4$, so
$$
(0,0,(1,2))^\perp\cong U\oplus U(2)\oplus <-4>\cong T_{ns1}\qquad
(\subset U\oplus U(2)\oplus A_2(-2)=T_{gen}).
$$

\subsection{Eckardt points}
As the Eckardt points of a cubic surface $S$ are the singular points of
its Hessian $Y$ which are on $S$ (\cite{S}, p.146), it is easy to find the
possible configurations of Eckardt points on a surface of type $ns1$.
We determine the transcendental lattices  of these surfaces below.

\subsection{Proposition}\label{nse}
Let $S_{ns1}$ be a cubic surface without a Sylvester form as in section
\ref{ns1}, from which we also adopt the notation and conventions:
$$
S_{ns1}:\qquad a_0x_0^3+x_1^3+x_2^3+x_3^3-3x_0^2(a_1x_1+a_2x_2+a_3x_3)=0.
$$
\begin{enumerate}
\item {}\label{a0=0}
In case $a_0=0$ and all $a_i^3$'s are distinct, the point $p_0$ is
the unique Eckardt point on $S_{ns1}$.
The general Hessian K3 surface has Picard number $18$
and transcendental lattice
$$
U(2)^2.
$$
\item{}\label{a0not0}
In case $a_0\neq 0$ but $a_i^3=a_j^3\neq a_k^3$ the point $q_k$ is
the unique Eckardt point on $S_{ns1}$. The general Hessian K3 surface has
Picard number $18$
and transcendental lattice
$$
U\oplus <-4>\oplus <4>.
$$
\item{} In case $a_0=0$ and $a_i^3=a_j^3\neq a_k^3$, the points $p_0,q_k$
are the only Eckardt points on $S_{ns1}$. The general Hessian K3 surface has
Picard number $19$
and transcendental lattice
$$
U(2)\oplus <4>.
$$
\item{}\label{ps}
In case $a_0\neq 0$ and $a_i^3=a_j^3= a_k^3$, the points $q_1,q_2,q_3$
are the only Eckardt points on $S_{ns1}$. The general Hessian K3 surface has
Picard number $19$ and transcendental lattice
$$
U\oplus  <12>.
$$
\item{}
In case $a_0=0$ and $a_i^3=a_j^3= a_k^3$, the points $p_0,q_1,q_2,q_3$
are the only Eckardt points on $S_{ns1}$. The general Hessian K3 surface has
Picard number $20$ and transcendental lattice
$$
A_2(2)=\mbox{$\left(
\begin{array}{cc}4&-2\\-2&4
\end{array}\right)$}.
$$
\end{enumerate}

\ts
We start with the case $a_0=0$ and $a_i^3=a_j^3= a_k^3$,
we denote the cubic surface by $S$.
(Note that the transformation $x_j\mapsto \omega^{n_j} x_j$, with $j=1,2,3$ and
with a cube root of unity $\omega$,
also gives an equation of type $ns1$ but with
$a_j\mapsto \omega^{-n_j}a_j$).
For every Eckardt point $p$ the
tangent plane to the cubic surface in $p$ intersects the Hessian
in a double line and two new lines. In particular, $p_0$ gives the pair of lines
lines $d^\pm$ on $S$ defined by
$\sum_{i=1}^3 x_i=0,\sum_{1\leq i< j\leq 3}x_ix_j=0$, and similarly
the $q_i$ give the pairs of lines $x_k+x_l=x_0^2+x_i^2=0$.
The lattice $L$ generated by a basis of $NS(X_{ns1})$, $d^+$ and the lines
$c^+_1:x_2+x_3=x_0+ix_1=0$, $c^-_2:x_1+x_3=x_0-ix_2=0$
has rank $20$ and discriminant $12$.
The only rank two, even, positive definite lattice with discriminant $3$
is $A_2$, but this cannot be the transcendental lattice of a Hessian by
Lemma \ref{slh}. Hence $L$ must be the N\'eron-Severi lattice of the Hessian K3
surface $X$ of $S$ and a computation of the discriminant group shows that
the transcendental lattice is $A_2(2)$. It is also not hard to find two
orthogonal $E_8$'s in the N\'eron-Severi group of $X$ and thus to find an
embedding of $NS(X)$ into the K3 lattice.

For the other cases it is easy to find the rank $r$ of the
transcendental lattice
from a count of the moduli and to find curves on the surface which span
a lattice of rank $22-r$. Then one can verify that these curves all lie
in the following sublattices of $L=NS(X)$, and using the embedding of
$NS(X)$ into the K3 lattice, one finds the transcendental lattices.
\begin{enumerate}
\item {} The N\'eron-Severi group is
$<(c^+_1-c^-_1),(c^+_2-c^-_2)-(c^+_3-c^-_3)>^\perp$
in $NS(X)$.
\item{} The N\'eron-Severi group is
$<(c^+_i-c^-_i)-(c^+_j-c^-_j),d^+-d^->^\perp$ in $NS(X)$.
\item{} The N\'eron-Severi group is $<(c^+_i-c^-_i)-(c^+_j-c^-_j)>^\perp$
in $NS(X)$.
\item{} The N\'eron-Severi group is $<d^+-d^->^\perp$ in $NS(X)$.
\end{enumerate}

\qed

\subsection{Remark}
The one dimensional family of Hessians in Proposition \ref{nse}.\ref{ps}
with transcendental lattice $U\oplus <12>$ was studied by Peters and Stienstra
in \cite{Stienstra}. They also observed the `cube' formed by 20 of the curves
in the N\'eron-Severi lattice (see the proof of Proposition \ref{tnsk}).
The total space of this family is a Calabi-Yau threefold (the Fermi threefold)
and was studied by Verrill \cite{Verrill}, section 4.

\section{The moduli space of cubic surfaces}\label{moduli}

\subsection{} Classical invariant theory shows that the moduli space
of cubic surfaces $\cM$ is isomorphic to the weighted projective space
$\PP(1,2,3,4,5)$.
Each point in $\cM$ corresponds to the isomorphism class
of a cubic surface with at most 4 nodes, except for one point.
That point  corresponds to the semi stable, non-stable, cubic surfaces,
the unique closed orbit in this set is the orbit of $t^3=xyz$.
In Theorem \ref{nsinv} we show that it maps to
$(8:1:0:0:0)\in \cM$
where we use Salmon's generating invariants (cf.\ section \ref{geninv}).

We recall the classical description of this moduli space using the
Sylvester forms. We discuss some divisors, and their classes
in the Chow group of the moduli space. We
determine the divisor
parametrizing cubic surfaces without a Sylvester form in Theorem \ref{nsinv}.
Finally we make some comments on
the singular locus of the moduli space.

\subsection{Invariants}
The ring of invariants of the action of $SL(4,\CC)$ on the space
of cubic forms in 4 variables is generated by the invariant polynomials
$I_n$ of degree $n$ for $n=8,16,24,32,40,100$.
Since $100$ is not divisible by $8$ and $I_{100}^2$ is polynomial in the
other generators, the moduli space of cubic surfaces (the $Proj$ of
the ring of invariants) is the weighted projective space $\PP(1,2,3,4,5)$.
Note that $I_8$ and $I_{100}$ are, up to a scalar multiple, unique whereas
$I_{16},\ldots,I_{40}$ are only unique up to the addition of
homogeneous weighted
polynomials of lower degree (for example, for any $a,b,c\in\CC$ $a\neq 0$,
the invariant $aI_{24}+bI_8^3+cI_8I_{16}$ can also be used as a generator
of degree $24$). We will use Salmon's convention to choose the generators.

\subsection{Generating invariants and the Sylvester form}\label{geninv}
The generating invariants are easily computed for a Sylvester form
$$
\sum_{i=0}^4 x_i=0,\qquad \sum_{i=0}^4\lambda_ix_i^3=0.
$$
Let $\sigma_i$ be the $i$-th symmetric function in
$\lambda_0,\ldots,\lambda_4$. Then (\cite{Salmon} p.\ 197):
$$
I_8=
\sigma_4^2-4\sigma_3\sigma_5,\quad
I_{16}=\sigma_5^3\sigma_1,\quad
I_{24}=\sigma_5^4\sigma_4,\quad
I_{32}=\sigma_5^6\sigma_2,\quad
I_{40}=\sigma_5^8.
$$

The quotient of $\PP^4_\lambda$, the parameter space of the Sylvester forms,
by the group $S_5$ is also a weighted projective space $\PP(1,2,3,4,5)_\sigma$
with coordinates the elementary symmetric functions. The formulas above define
a birational isomorphism:
$$
\phi:\PP(1,2,3,4,5)_\sigma --\rightarrow\PP(1,2,3,4,5)_I,
\qquad (\sigma_1:\ldots:\sigma_5)\longmapsto (I_8:\ldots:I_{40}).
$$
The base locus of $\phi$ is given by $\sigma_4=\sigma_5=0$.
The birational inverse of $\phi$ is:
$$
\psi:\PP(1,2,3,4,5)_I --\rightarrow\PP(1,2,3,4,5)_\sigma,
$$
$$
(I_8:\ldots:I_{40})\longmapsto
(I_{16}:I_{32}:(I^2_{24}-I_8I_{40})/4:I_{24}I_{40}:I_{40}^2),
$$
the base locus of $\psi$ is the point $(1:0:0:0:0)$ which corresponds to
the Fermat cubic surface, see Theorem \ref{nsinv}.

\subsection{The boundary divisor}\label{boundary}
The locus of singular cubic surfaces, the boundary of the moduli space,
is the divisor in $\PP(1,2,3,4,5)$ defined by:
$$
(I_8^2-2^6I_{16})^2=2^{14}(I_{32}+2^{-3}I_8I_{24}).
$$
The degree of this invariant is 32
(cf.\ \cite{S}, App.\ III p.\ 179, \cite{GKZ}, p.\ 38), the formula
in \cite{Salmon}, p.\ 198,  omits the exponent $-3$, but we verified that
upon substituting the $\sigma_i$ in the $I_n$ in our formula
one finds the expression given in section \ref{discriminant}.
As $I_{40}$ does not appear in this formula, the boundary divisor is a
cone over a $\PP(1,2,3,4)$.

\subsection{The tritangent divisor}\label{tritangent}
The tritangent divisor parametrizes cubic surfaces with an Eckardt point.
The general such surface has a Sylvester form in which two of the $\lambda_i$
coincide. Thus this divisor is defined by the discriminant of the
polynomial $\prod_{i=0}^4 (x-\lambda_i)$, which is the invariant $I_{100}$
(\cite{Salmon}, p.\ 197).
The discriminant is a polynomial of weight 20 in the
symmetric functions in the $\lambda_i$ and gives, upon pull-back along the map
$\psi$ from section \ref{geninv}, a polynomial of weight $320=20\cdot 16$
in the basic invariants $I_8,\ldots,I_{40}$.
This polynomial has a factor $I_{40}^3$, the other
factor is $I_{100}^2=f(I_8,\ldots,I_{40})$ which defines
the tritangent divisor in $\PP(1,2,3,4,5)$.

\subsection{Theorem}\label{nsinv}
The divisor in the moduli space of cubic surfaces $\cM\cong\PP(1,2,3,4,5)_I$,
which parametrizes cubic surfaces which do not admit a Sylvester form
is defined by $I_{40}=0$.

\noindent
The subvariety parametrizing the surfaces of type $ns2$ is defined by
$I_{24}=I_{40}=0$.

\noindent
The subvariety parametrizing the cyclic surfaces is defined by
$I_{24}=I_{32}=I_{40}=0$.

\noindent
The point $I_{16}=I_{24}=I_{32}=I_{40}=0$ parametrizes the
Fermat cubic surface.

\noindent
The non-stable, semi-stable cubic surface defined by $t^3=xyz$
maps to the point $(8:1:0:0:0)\in\cM$.

\ts
A family of five linear forms where $x_0,\ldots,x_3$
are coordinates on a $\PP^3$
and where two of the linear forms coincide in the limit $t\rightarrow 0$ is:
$$
x_0(t)=x_0,\quad x_i(t)=a_itx_i\quad(i=2,3,4),\qquad
x_4(t)=-x_0-t(a_1x_1+a_2x_2+a_3x_3).
$$
Note that any four of the five linear forms are independent when the
$a_i\neq 0$ and $t\neq 0$ and that $\sum x_i(t)=0$.
The family of cubic surfaces, defined by Sylvester forms:
$$
S_t:\quad\left\{\begin{array}{ll}
(a_0+t^{-1})x_0(t)^3+(a_1t)^{-3}x_1(t)^3+(a_2t)^{-3}x_2(t)^3+
(a_3t)^{-3}x_3(t)^3+t^{-1}x_4(t)^3=0,\\
\sum_{i=0}^4x_i(t)=0\end{array}\right.
$$
has as limit the cubic surface $S_{ns1}$ defined by:
$$
S_{ns1}:=S_0:\qquad a_0x_0^3+x_1^3+x_2^3+x_3^3-3x_0^2(a_1x_1+a_2x_2+a_3x_3)=0.
$$
Replacing $x_0$ by $-x_0$ we get the cubic form:
$$
F=x_1^3+x_2^3+x_3^3-x_0^2(a_0x_0+3a_1x_1+3a_2x_2+3a_3x_3)
$$
considered in \ref{ns1}.

To find the invariants of the general surface of type $ns1$
we compute the invariants of the surface $S_t$
in the one parameter family $S_t$
and take the limit $t\rightarrow 0$. The result is:
$$
[S_{ns1}]=[\lim_{t\rightarrow 0} S_t]
=(-4\rho_1+a_0^2:\rho_2:2\rho_3:\rho_1\rho_3:0)
\qquad (\in\PP(1,2,3,4,5)_I),
$$
where $\rho_i$ is the $i$-th elementary symmetric function in
$a_1^3,a_2^3,a_3^3$.
From the equation of $S_{ns1}$ it is clear that its isomorphism
class depends only on the symmetric functions of $a_1^3,a_2^3,a_3^3$.
Thus the three dimensional family $ns1$ maps to the divisor $I_{40}=0$
which shows that this divisor is the closure of the locus of surfaces
of type $ns1$. Note that for a general point of $I_{40}=0$, in particular
one for which the other coordinates are non-zero, one can recover $\rho_2$ and
$\rho_3$ from the second and third coordinate, the fourth coordinate
then determines $\rho_1$ and
finally $a_0$ is recovered from the first coordinate.

A one parameter family where three of the linear forms coincide
in the limit $t\rightarrow 0$ is:
$$
x_0(t)=-x_3-t^2x_0,\quad x_1(t)=\mu t^2x_1,\quad x_2=t^2x_2,\quad
x_3=2x_3+tx_0
$$
and $x_4(t)=-(x_0(t)+\ldots+x_3(t))$. Again any four of the 5 $x_i(t)$ are
independent for $t\neq 0$.
The family of cubic surfaces, defined by Sylvester forms:
$$
S_t:\quad\left\{\begin{array}{ll}
t^{-2}x_0(t)^3+(\mu t^2)^{-3}x_1(t)^3+t^{-6}x_2(t)^3+
(\mbox{$\frac{\lambda}{4}+\frac{1}{4t^2}$})x_3(t)^3+t^{-2}x_4(t)^3=0,\\
\sum_{i=0}^4x_i(t)=0\end{array}\right.
$$
has as limit the cubic surface $S_{ns2}$ defined by (after substituting
$\mbox{$\frac{1}{\sqrt{2}}$}x_0$ for $x_0$):
$$
S_{ns2}:=S_0:\qquad G=x_1^3+x_2^3+2\lambda x_3^3-3x_3(\mu x_1x_3+x_2x_3+x_0^2)=0.
$$
The same procedure as above gives:
$$
[S_{ns2}]
=(-8\lambda:1+\mu^3:0:\mu^3:0)
\qquad (\in\PP(1,2,3,4,5)_I),
$$
where $\lambda,\mu$ are the coefficients of the equation of $S_{ns2}$.

Finally we consider the cyclic surfaces. These can be obtained from
$S_{ns1}$ with one of the $a_i=0$, but the following family has a particularly
nice limit, a sum of cubes of linear forms.
Consider the one parameter family of Sylvester forms:
$$
x_i(t)=(1+t)x_i,\quad(i=0,1,2),\qquad
x_4(t)=tx_4,\qquad x_3(t)=-(x_0(t)+x_1(t)+x_2(t)+x_4(t)).
$$
Any four of the five $x_i(t)$ are
independent for $t\neq 0$.
The family of cubic surfaces, defined by Sylvester forms:
$$
S_t:\quad
\lambda_0x_0(t)^3+\ldots +\lambda_3 x_3(t)^3+
(\mbox{$\frac{\lambda_4}{t^3}$})x_4(t)^3=0,    \quad
\sum_{i=0}^4x_i(t)=0,
$$
has as limit the cyclic cubic surface $S_{cyc}$ defined by
$$
S_{cyc}:=S_0:\qquad
G=\lambda_4x_4^3
-\lambda_3(x_0+x_1+x_2)^3+
\lambda_0x_0^3+\lambda_1x_1^3+\lambda_2x_2^3=0.
$$
Computing the invariants as before gives:
$$
[S_{cyc}]=(\tau_3^2-4\tau_2\tau_4:\tau_4^3:0:0:0)\qquad (\in\PP(1,2,3,4,5)_I),
$$
where $\tau_i$ is the $i$-th elementary symmetric function in
$\lambda_0,\ldots,\lambda_3$ (note that $\lambda_4$ does not appear
as is obvious from the equation of $S_{cyc}$).
The Fermat cubic is the surface $S_{cyc}$ with
$\lambda_3=0$, hence $\tau_4=\prod_{i=0}^3\lambda_i=0$.

Putting $\lambda_i=1$, $0\leq i\leq 3$, we obtain the cubic surface
$x_4^3=3 (x_0 + x_1) (x_0 + x_2) (x_1 + x_2)$, which is isomorphic
to the strictly semi-stable surface $t^3=xyz$, and it defines the point
$(-8:1:0:0:0)=(8:1:0:0:0)\in\cM=\PP(1,2,3,4,5)$.
Note that this is the unique point in the intersection of the boundary
divisor, defined by the equation from \ref{boundary},
with the curve $I_{24}=I_{32}=I_{40}=0$ parametrizing cyclic surfaces.
\qed

\subsection{The Kummer divisor}
The Kummer K3 surface of a principally polarized abelian surface is
the Hessian K3 surface of a cubic surface.
The divisor in $\cM$ which we obtain in this way will be called the
Kummer divisor.
Rosenberg (\cite{Rosenberg}, Cor.\ 1.2)
proved that the Kummer divisor in $\cM$ is defined by:
$$
I_8I_{24}+8I_{32}.
$$
The N\'eron-Severi group of the Kummer surface has rank $17$ and has
transcendental lattice isomorphic to $T_{kum}=U(2)\oplus U(2)\oplus <-4>$.
The lattice $T_{kum}$ is the orthogonal complement of an element
$t$ with $t^2=-12$ in $T_{gen}$.

\subsection{The locus $I_{40}=I_{100}^2=0$}
The intersection of the tritangent divisor
$I_{100}^2=0$
and the non-Sylvester divisor $I_{40}=0$ consists of three irreducible
components (of dimension two).
In fact, putting $I_{40}=0$ in $I_{100}^2=f(I_8,\ldots,I_{40})$
one finds (up to scalar multiple):
$$
f(I_8,\ldots,I_{32},0)=I_{24}^3(I_8I_{24}+8I_{32})g(I_8,\ldots,I_{32})
$$
where $g$ is given by:
$$
g=16I_{16}^3I_{24}^2 + 27I_{24}^4 - 72I_{16}I_{24}^2I_{32} -
16I_{16}^2I_{32}^2 +   64I_{32}^3.
$$

The first factor, $I_{24}=0$, corresponds to the non-Sylvester surfaces
of type $ns2$, cf. Theorem \ref{nsinv}.
These surfaces do indeed have an Eckardt point:
$p_0=(1:0:0:0)$
(notation as in section \ref{ns2}) is singular on the Hessian and lies
on $S_{ns2}$, hence it is an Eckardt point.

The second factor is the intersection of $I_{40}=0$ with the Kummer divisor.
Parametrizing the non-Sylvester divisor as in the proof of Theorem \ref{nsinv},
this component is given by $a_0=0$ and these surfaces in fact have an Eckardt
point (cf.\ Proposition \ref{nse}.\ref{a0=0}).

The last factor pulls-back to the discriminant
of the polynomial $(x-a_1^3)(x-a_2^3)(x-a_3^3)$, hence gives the surfaces
with Eckardt point considered in Proposition \ref{nse}.\ref{a0not0}.

\subsection{The singular points of the moduli space}
The singular locus of a weighted projective space was determined in
\cite{dimca}, Prop.\ 7. For $\cM=\PP(1,2,3,4,5)$ we get:
$$
\PP(1,2,3,4,5)_{sing}=\{(0:0:1:0:0)\}\cup\{0:0:0:0:1)\}\cup
\{(0:a:0:b:0)\},
$$
with $a,b\in\CC$ not both zero.

Let $\PP(q_1,\ldots,q_n)$ be
a weighted projective space
where we (may) assume that any $n-1$ of the $q_i$ are relatively prime.
This algebraic variety is covered by affine
open subsets $U_i$ which are cyclic quotients of vector spaces
(\cite{dolwps}, Proposition 1.3.3 and its proof):
$$
U_i=\{(x_1:\ldots:x_n)\in\PP(q_1,\ldots,q_n):x_i\neq 0\}\;\cong\;
\CC^{n-1}/\mu_{q_i},
$$
where the action of a $q_i$-th root of unity $\zeta$ on $\CC^{n-1}$
is given by
$$
(x_1,\ldots,x_{i-1},x_{i+1},\ldots,x_n)\longmapsto
(\zeta^{q_1}x_1,\ldots,\zeta^{q_{i-1}}x_{i-1},
\zeta^{q_{i+1}}x_{i+1},\ldots,\zeta^{q_n}x_n).
$$
In particular, $U_1\subset\PP(1,2,3,4,5)$ is isomorphic to
$\CC^4$, $U_3,U_5$ are the quotient of $\CC^4$
by an automorphisms of order three and five respectively
and both have an isolated singularity.

\subsection{The cubic surfaces in the singular locus.}
The singular locus of $\cM$ parametrizes, as expected, certain cubic surfaces
with non-trivial automorphism groups.

The cubic surface corresponding to $(0:0:1:0:0)$ has $I_{40}=0, I_{24}\neq 0$,
hence it is of
non-Sylvester type $ns1$ (Proposition \ref{nsinv}), moreover, it is a
singular point (section \ref{boundary}). The parametrization given in
the proof of Theorem \ref{nsinv} shows that $\rho_1=\rho_2=a_0=0$,
from which one finds that $a_i^3=\omega^ia_3^3$ for a primitive cube of
unity $\omega$ and the equation of the surface is:
$
x_1^3+x_2^3+x_3^3-3x_0^2(x_1+\epsilon x_2+\epsilon^2x_3)=0
$
where $\epsilon$ is a primitive $9$th root of unity. A change of variables
gives the equation ($\omega=\epsilon^6$):
$$
x_1^3+\omega x_2^3+\omega^2 x_3^3-3x_0^2(x_1+x_2+x_3)=0.
$$
This surface has two singular points, $p_\pm=(\pm 1:1:\omega:\omega^2)$.
It has an automorphism of order three given by
$$
(x_0:x_1:x_2:x_3)\longmapsto (\omega x_0:x_2:x_3:x_1).
$$

The cubic surface corresponding to $(0:0:0:0:1)$ is singular and has
the Sylvester form:
$$
\sum_{i=0}^4 \eta^ix_i^3=0,\qquad \sum_{i=0}^4 x_i=0,
$$
where $\eta$ is a primitive $5$-th root of unity.
This surface has the unique singular point $(1:\eta^2:\eta^4:\eta:\eta^3)$,
it is an ordinary double point.
The map:
$$
(x_0:\ldots:x_3:x_4)\longmapsto (x_1:\ldots:x_4:x_0)
$$
is an automorphism of order five of this surface.

The rational curve defined by $I_8=I_{24}=I_{40}=0$ parametrizes
cubic surfaces of non-Sylvester type $ns2$ (Theorem \ref{nsinv}):
$$
S_\mu:\quad x_1^3+x_2^3-3x_3(\mu x_1x_3+x_2x_3+x_0^2)=0
$$
and $\mu^3/(1+\mu^3)^2=b/a^2$ (or, if $a=0$, $\mu^3=-1$).
The point $(1:0:0:0)$ is an Eckardt point.
Each of these surfaces has an automorphism of order $4$ given by:
$$
(x_0:x_1:x_2:x_3)\longmapsto (ix_0:x_1:x_2:-x_3).
$$
This family was studied for example in \cite{S}, p.\ 151,
\cite{Naruki}, Thm.\ 2, p.\ 3, \cite{TH}, Thm.\ 5.3. Any surface
in the family has the property that its automorphism group modulo
the subgroup generated by the involutions defined by the Eckardt points
is non-trivial. This quotient group has two elements except if $\mu^6=1$.
When $\mu=\omega$, a primitive cube root of unity,
the surface has two singular points
$(1:0:\omega^2:\pm 1)$ and it has another Eckardt point: $(0:1:-\omega:0)$.
In case $\mu^3=-1$ one obtains a surface
with quotient group isomorphic to $\ZZ/4\ZZ$. This surface defines the point
$(0:0:0:1:0)\in\cM$.

The surface with $\mu=0$, which corresponds to $(0:1:0:0:0)\in\cM$,
has affine equation (putting $x_0=y,x_1=t,x_2=x,x_3=1$ in the equation above):
$$
S_0:\quad t^3+3y^2+x^3-3x=0
$$
which shows that it is a cyclic surface (as expected, see Theorem \ref{nsinv})
which branches over the smooth
elliptic curve $E:\;3y^2+x^3-3x=0$ in $\PP^2$, it has $j(E)=1728$,
in fact the automorphism of $S_0$ of order 4 above induces one on $E$ with
a fixed point.

\end{document}